\documentclass[12pt]{article}
\usepackage[latin1]{inputenc} 
\usepackage{latexsym}
\usepackage{csvsimple}
\usepackage{bbm}              
\usepackage{dsfont}           
\usepackage{graphicx}         
\usepackage{epsfig}           
\usepackage{amsmath}
\usepackage{amsthm}         
\usepackage{amsfonts}
\usepackage{amssymb}
\usepackage{empheq}
\usepackage{color}
\usepackage{array}
\usepackage{amsmath}
\usepackage{framed}
\usepackage{fancyhdr}
\usepackage{array,multirow,makecell}
\usepackage{multicol}
\usepackage{titlesec}
\usepackage{tikz}
\usepackage{times}
\usepackage{graphics}
\usepackage[T1]{fontenc}
\usepackage{algorithm}
 \usepackage{algorithmic}
 \usepackage{caption}
\usepackage[colorlinks,citecolor=red,urlcolor=blue,bookmarks=false,hypertexnames=true]{hyperref}
\newcolumntype{L}[1]{>{\raggedleft\arraybackslash   }b{#1}}

\theoremstyle{plain}
\newtheorem{theorem}{Theorem}[section]

\theoremstyle{definition}
\newtheorem{definition}[theorem]{Definition}

\begin{document}


%

\begin{center} {\Large A Sequential Descent Method for Global Optimization}
\end{center}
\begin{center}
{\footnotesize Mohamed Tifroute$^1$,  Anouar Lahmdani$^2$ and Hassane Bouzahir$^3$ \\
E-mails: mohamed.tifroute@gmail.com, anouar.lahmdani@gmail.com
}\end{center}

\vspace{30pt}

\begin{center} {\footnote{Higher School of Technology - Guelmim , Ibn Zohr University, Morocco}, \footnote{
FSA - Ait Melloul, Ibn Zohr University, Morocco} and \footnote{ENSA, Ibn Zohr University, Morocco } }
\end{center}

\clearpage

\begin{abstract}
In this paper, a sequential search method for finding the global minimum of an objective
function is presented, The descent gradient search is repeated  until the global minimum is obtained.  The global minimum is located by a process of finding progressively better local minima. We determine the set of points of intersection between the curve of the function  and the horizontal plane which contains the local minima previously found. Then, a point in this set with the greatest descent slope is chosen to be a initial point for a new descent gradient search. The method has the descent property and the convergence is monotonic. To demonstrate the effectiveness of the proposed sequential descent method, several non-convex multidimensional optimization problems are solved. Numerical examples show that the global minimum can be sought by the proposed method of sequential descent.
\end{abstract}
\vspace{30pt}

\begin{sloppypar}
\section{Introduction}
Multi-dimensional non-convex continuous optimization problems are important in
many practical applications.   Many approaches which are supported by relevant
convergence analysis are for finding local minima only (Rubinov \cite{1.}). However,
many local minima are useless in practice, as their corresponding cost values are
too much inferior to the global minimum cost value. There are several proposed stochastic
optimization methods  to solve for the global minimum. Global optimization is a NP complete problem and heuristic approaches like Genetic Algorithms (GAs) and Simulated Annealing
(SA) have been used historically to find near optimal solutions \cite{1., 2.}.
The Particle Swarm Optimization (PSO) algorithm proposed by Kennedy and Eberhart in 1995 \cite{3., 5.} is a metaheuristic algorithm based on the concept of swarm intelligence capable of solving complex mathematics problems existing in engineering. It has been shown that the PSO algorithm is easy to implement  and it  converges
faster than traditional techniques like GAs for a wide variety of
benchmark optimization problems \cite{6., 7.}. Heuristics method are very expensive to apply. Therefore, methods which are hybriding different type of algorithms are becoming more popular.
One method is to use gradient-type procedures coupled with certain auxiliary
functions to move successively from one local minimum to another better one.
This includes the tunnelling method. For more
details, see,Levy an Montalvo \cite{10.}, Yao \cite{20.} and Liu and Teo, \cite{12.}.
 Yiu et al. \cite{14.} propose a hybrid algorithm in which SA and the Descent Method (DM) are used. SA is used in the global search phase, and the DM is used in the local search phase. In \cite{15.},
metaheuristic methods exploiting a chaotic system based on the DM have been presented.

In general, an efficient global optimization algorithm would have the capacity
of overcoming the problem of local minima, and assuring the speed of convergence to approach stationary
points.  On the part of continuous decision variables, the stochastic optimization
approach gives a good methodology to get away from stationary points, but
it is computationally intensive to be practicable. 	Among the  argumentations for this problem, the fact that
the method is very slow when it strives to nearing or descents to stationary points.
However, an analytic approach based on the gradient information is much
more efficient in finding a stationary point. In this paper, a sequential technique
is proposed consisting in repeating the classical method of the descent gradient until the global minimum is obtained.  The global minimum is located by a process of finding progressively better local minima. We determine the set of points of intersection between the curve of the function and the horizontal plane which contains the previously found local minima .
Accordingly, we propose to repeat this process and choose a descent point from the intersection between the horizontal plan containing the previously converged local minima and the curve of the function. The probability of finding a better descent point is much larger then finding a better local minimum. Thus, it is much more efficient computationally.
The advantage of the proposed sequential descent method is
that the convergence is monotonic. The decrease in the objective function after executing
each descent gradient search might be very small. But, it is sufficient to detour previously converged local solutions.
To demonstrate the effectiveness of the proposed sequential method, several multidimensional non-convex optimization problems are solved. For each example, the proposed sequential descent method locates the corresponding global solution.
\section {Deterministic Optimization Algorithm}
Deterministic optimization algorithms are generally based on the gradient of the objective function with respect to design variables.
The deterministic method, when applied to non-linear non-convex minimization
problems, usually requires implementing an iterative method,
which, will optimistically converge to a local minimum of the objective function, after a certain number of iterations. The iterative method can be written as follows.
\begin{equation}\label{22}
x^{k+1} = x^k + \lambda_k d^k
\end{equation}
where $x \in X \subset \mathbb{R}^p$ is the vector of design variables, $\lambda$ is the
search step size, $d$ is the direction of descent, $k$ is the number
of iterations and $X$ is the search space.
An iteration step is acceptable if $f(x^{k+1}) < f(x^k)$. The direction of descent $d$ will generate an acceptable step if and only if there exists a positive definite matrix $\textbf{M}$, such that
$d = -\textbf{M}\nabla f$, where $\nabla f$ is the gradient of $f$.
Such requirement results in directions of descent that
form an angle greater than $90^{°}$ with the gradient direction.
A minimization method in which the directions are obtained
as above is called an acceptable gradient method.\\
\textbf{Steepest Descent Method} For the steepest descent
method, the direction of descent is given by
\begin{equation}\label{33}
d^k = -\nabla f(x^k)
\end{equation}
\textbf{Line Search} For a descent direction $d^k$, the line search procedure determines the solution of the following one-dimensional optimization problem
\begin{equation}\label{44}
\lambda_k = \displaystyle { argmin_{\lambda}} f(x^k + \lambda d^k)
\end{equation}
and takes $\lambda_k$ as a step size. In fact, there are several
approaches that define conditions for a step size that could
be obtained as an approximate solution of \ref{55}. One of them
is the Armijo rule,
\begin{equation}\label{55}
f(x^k + \lambda_k d^k) \leq f(x^k) + c_1 \lambda_k \nabla f(x^k)^T d^k
\end{equation}
where $c_1$ is a small positive constant.
While the line search procedure works well with algorithms of the first order (i.e., algorithms using function and gradient values only) the convergence rate is at most linear for deterministic rate.
\begin{definition}
Let  a vector $v \in \mathbb{R}^n$ be given by $v=(v_1,v_2,\cdots,v_n)$, $v$ is said
to be strictly negative if all component $v_1,v_2,\cdots,v_n$ are strictly negative.
\end{definition}
\begin{algorithm}
 \caption{SGD (Sequential Gradient  Descent)}
 \begin{algorithmic}[1]

 \STATE Generate $x^{(0)}$ randomly and evaluate $f(x^{(0)})$. Set $k = 0$;\\

 \STATE  Solve for the local minimum of $f(x)$ using a gradient-based minimization method
with $x^{(k)}$ as the initial guess to give $x^{(k)*}$ such that $|f(x^{(k)*}) - f(x^{(k)})| \leq \epsilon_k$,
where $\epsilon_k$ is a positive parameter;\\
 \STATE Define $(P_k)$ the horizontal plan containing $f(x^{(k)*})$;\\
 \STATE Find the set $I_k= (P_k)\bigcap (Cf)$;\\
 \textbf{if} for all $g_k \in I_k$, $\nabla f(g_k)=0$, return $x^{(k)*}$;
\textbf{ else};\\
 \STATE  Put :$I_{k}^{-}=\{g_k \in I_k / g_k <0 \}$, choose $x^{(k+1)}= \displaystyle argmax_{g_k \in I_{k}^{-}}\{\|\nabla f(g_k)\| \} $;\\
let $k := k+1$\\
go to Step 2.
 \end{algorithmic}
 \end{algorithm}
 \newpage
 \section{Numerical Examples }
\begin{figure}
  \centering
  \includegraphics[width=12cm]{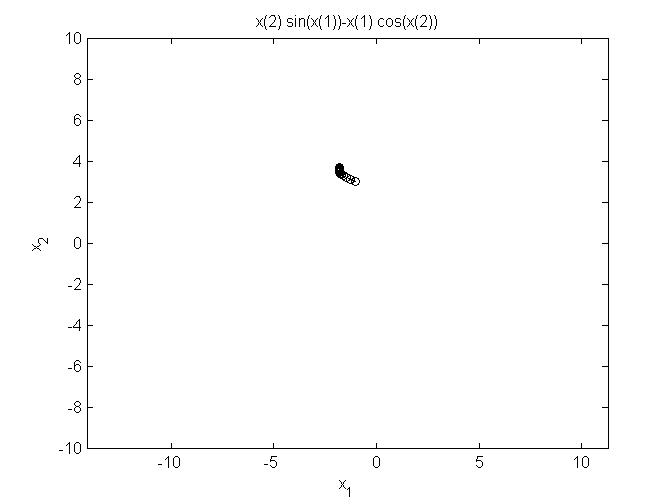}\\
  \caption{First local optimal solution $x^{1,*}$ }\label{local1}
\end{figure}

\begin{figure}
  \centering
  \includegraphics[width=12cm]{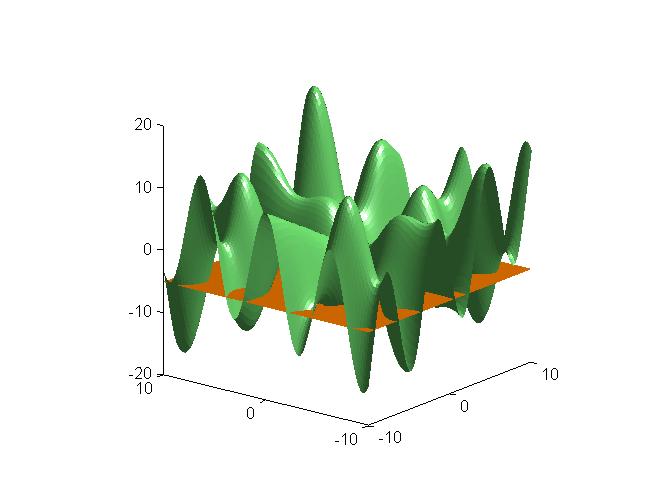}\\
  \caption{Intersection between the curve of the function and the horizontal plan containing $f(x^{1,*})$ }\label{zzz}
\end{figure}

\begin{figure}
  \centering
  \includegraphics[width=12cm]{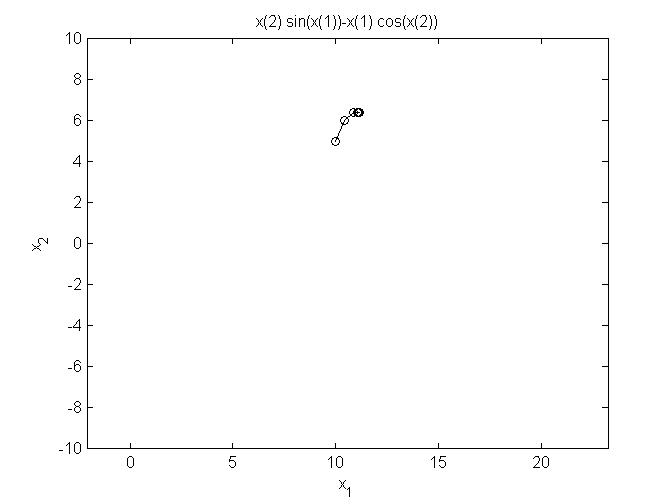}\\
  \caption{Second local optimal solution $x^{2,*}$}\label{local2}
\end{figure}

\begin{figure}
  \centering
  \includegraphics[width=12cm]{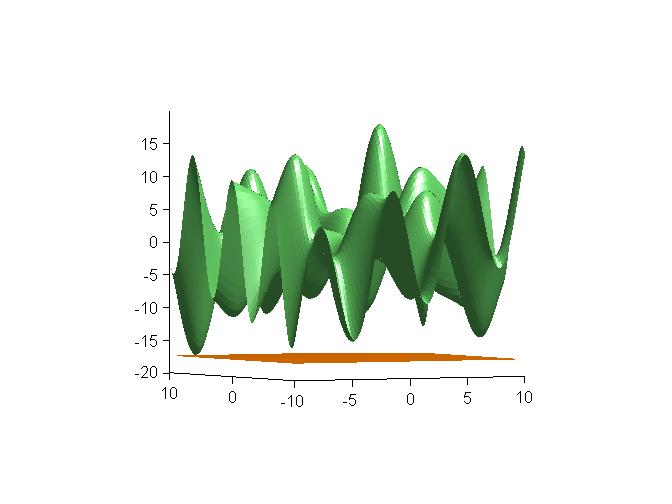}\\
  \caption{Intersection between the curve of the function and the horizontal plan containing $f(x^{2,*})$}\label{sss}
\end{figure}

{\bf Example 1:} \\
Consider the following problem:\\
Minimize $x_2 sin(x_1)-x_1 cos(x_2)$\\
$-5 \leq x_j \leq 5, j = 1, 2$.\\
\textbf{Step 1}: Initial start point $(-1, 3)$, local optimal solution is $(-1.7992, 3.7137)$ and $f(x^{1,*})=-5.1300$ (see figures \ref{local1}   and \ref{zzz}).\\
\textbf{Step 2}: Best new start point is : $(10 , 4.958375346)$, the new local optimal solution is $(11.1525 , 6.3719)$ and $f(x^{2,*})=-17.4022$ (see figure \ref{local2} ).\\
The intersection between the curve of the function and the horizontal plane which contains $f(x^{2,*})$  is reduced to a single point (see figure \ref{sss}), and as the gradient of the function at this point is close to zero, we deduce that the solution $x^{2,*}$  is the global minimum of the function.


\begin{figure}
  \centering
  \includegraphics[width=12cm]{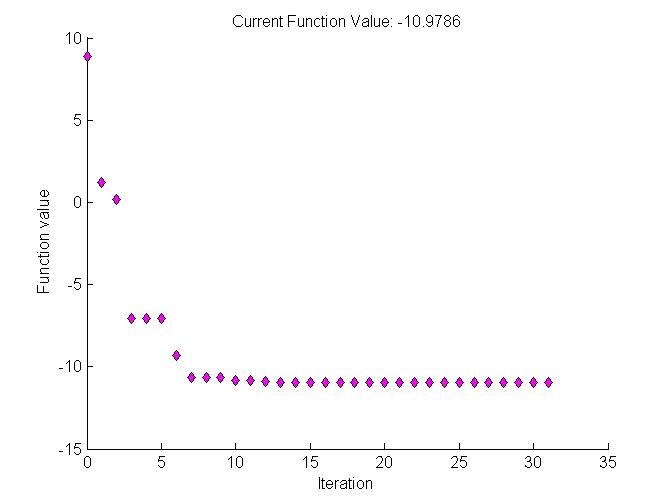}\\
  \caption{First local optimal solution $x^{1,*}$ }\label{local11}
\end{figure}

\begin{figure}
  \centering
  \includegraphics[width=12cm]{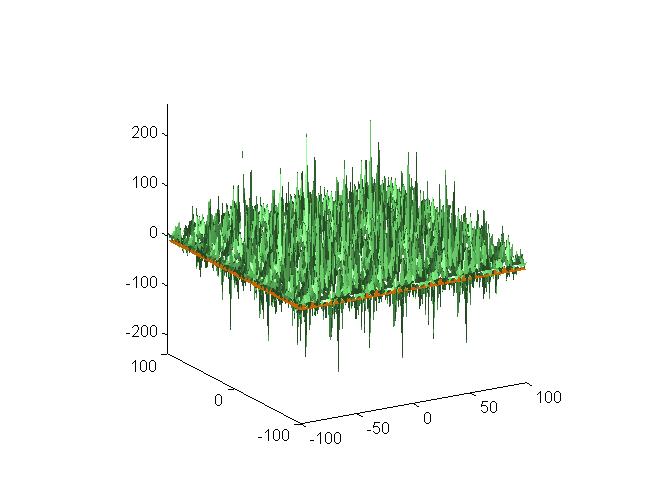}\\
  \caption{Intersection between the curve of the function and the horizontal plan containing $f(x^{1,*})$ }\label{zzz1}
\end{figure}

\begin{figure}
  \centering
  \includegraphics[width=12cm]{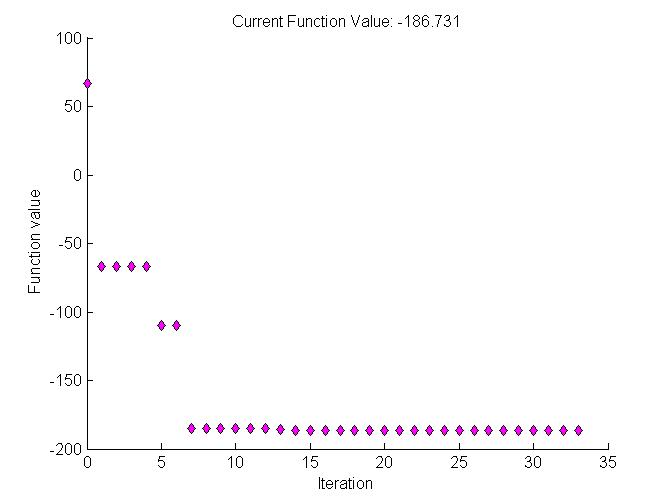}\\
  \caption{Second local optimal solution $x^{2,*}$}\label{local111}
\end{figure}

\begin{figure}
  \centering
  \includegraphics[width=12cm]{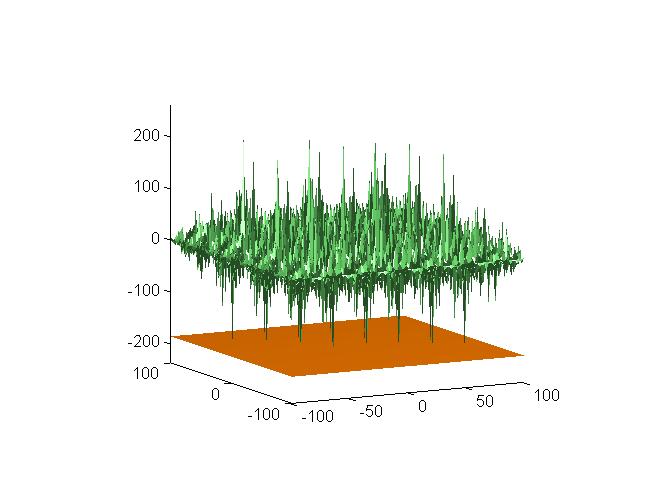}\\
  \caption{Intersection between the curve  of the  function and the horizontal plan containing $f(x^{2,*})$ }\label{bz}
\end{figure}

\newpage
\textbf{Example 2:}\\
The two-dimensional Shubert function (Shubert, 1972 )
$$
\begin{array}{l}
f\left(x_{1}, x_{2}\right)=\left(\sum_{i=1}^{5} i \cos \left[(i+1) x_{1}+i\right]\right)\left(\sum_{i=1}^{5} i \cos \left[(i+1) x_{2}+i\right]\right) \\
+\frac{1}{2}\left(\left(x_{1}+1.42513\right)^{2}+\left(x_{2}+0.80032\right)^{2}\right), \quad -100 \leqslant x_{i} \leqslant 100, \quad i=1,2
\end{array}
$$
\textbf{Step 1}: Initial start point $( 7 ,7)$, local optimal solution is $(7.603 , 7.105)$ and $f(x^{1,*})=-10.978558554610121$ (see figures \ref{local11} and \ref{zzz1}).\\
\textbf{Step 2}: Best new start point is  $(-13.13131313,	17.88967353)$, the new local optimal solution is $(-13.991  , 18.049)$ and $f(x^{2,*})=-186.73090665088228$ (see figure \ref{local111}).\\
The intersection between the curve of the function and the horizontal plane which contains $f(x^{2,*})$  is reduced to a finite number of points, and as the gradient of the function at this points is close to zero, we deduce that the solution $x^{2,*}$  is the global minimum of the function (see figure \ref{bz}).

In order to verify the quality of our algorithm, we compared our result to those of Yiu et al. \cite{14.} which obtains global minimum after 19 local search when our algorithm obtains the same result in only 2 sequential local search.


\newpage
{\bf Conclusion:} \\
In this paper, we have proposed a new sequential search method for finding the global minimum of an objective
function. Numerical results show that the proposed sequential method is promising for solving multidimensional non-convex continuous optimization problems.


%

 \end{sloppypar}

 \newpage

\bibliographystyle{plain}


\end{document}